\documentclass[12pt]{amsart}
\usepackage{ae,aecompl,amsbsy,amsfonts,amsmath,amsopn,amssymb,amsthm,amsxtra}
\usepackage{bm,bbm,blindtext,color,enumerate,etex,float,graphicx}
\usepackage{mathrsfs,multirow,pbox,pgf,stmaryrd,tikz,tikz-cd,upref,url,xcolor}
\usepackage[utf8]{inputenc}
\usepackage[T1]{fontenc}
\usepackage[english]{babel}
\usepackage[all,color]{xy}
\usepackage[colorlinks,final,backref=page,hyperindex]{hyperref}
\usepackage[square,sort,comma,numbers]{natbib}
\usetikzlibrary{shapes.multipart,positioning,decorations.pathreplacing,calc,intersections,folding}

\newtheorem{theorem}{Theorem}
\newtheorem{proposition}[theorem]{Proposition}
\newtheorem{lemma}[theorem]{Lemma}
\newtheorem{corollary}[theorem]{Corollary}
\newtheorem*{theorem*}{Theorem}

\theoremstyle{definition}

\theoremstyle{remark}
\newtheorem{remark}[theorem]{Remark}

\DeclareMathOperator{\rank}{rk}

\DeclareMathOperator{\supp}{s}

\newcommand{\xyinc}{\ar@{^{(}->}}
\newcommand{\xyrinc}{\ar@{_{(}->}}
\newcommand{\xyonto}{\ar@{->>}}
\newcommand{\xytwo}{\ar@{<->}}

\newcommand{\qand}{\quad\text{and}\quad}

\newcommand{\Kb}{\Bbbk}
\newcommand{\Rb}{\mathbb{R}}

\newcommand{\arr}{\mathcal A} %hyperplane arrangement 
\newcommand{\carr}{\mathcal C} %coordinate arrangement
\newcommand{\Hr}{\mathrm{H}} %hyperplane
\newcommand{\face}{\Sigma}
\renewcommand{\flat}{\Pi}
 
\newcommand{\minface}{O} %central face of an arrangement.

\newcommand{\maxflat}{\top}
\newcommand{\join}{\vee} %join 
\newcommand{\cpoly}{\chi} %characteristic polynomial

\newcommand{\rX}{\mathrm{X}} 
\newcommand{\rY}{\mathrm{Y}} 
 
\newcommand{\B}[1]{\mathtt{#1}} %The letter B stands for basis. 

\tikzset{
  anticlockwise arc centered at/.style={
    to path={
      let \p1=(\tikztostart), \p2=(\tikztotarget), \p3=(#1) in
      \pgfextra{
        \pgfmathsetmacro{\anglestart}{atan2(\y1-\y3,\x1-\x3)}
        \pgfmathsetmacro{\angletarget}{atan2(\y2-\y3,\x2-\x3)}
        \pgfmathsetmacro{\angletarget}%
        {\angletarget < \anglestart ? \angletarget+360 : \angletarget}
        \pgfmathsetmacro{\radius}{veclen(\y1-\y3,\x1-\x3)}
      }
      arc(\anglestart:\angletarget:\radius pt) -- (\tikztotarget)
    },
  },
  clockwise arc centered at/.style={
    to path={
      let \p1=(\tikztostart), \p2=(\tikztotarget), \p3=(#1) in
      \pgfextra{
        \pgfmathsetmacro{\anglestart}{atan2(\y1-\y3,\x1-\x3)}
        \pgfmathsetmacro{\angletarget}{atan2(\y2-\y3,\x2-\x3)}
        \pgfmathsetmacro{\angletarget}%
        {\angletarget > \anglestart ? \angletarget - 360 : \angletarget}
        \pgfmathsetmacro{\radius}{veclen(\y1-\y3,\x1-\x3)}
      }
      arc(\anglestart:\angletarget:\radius pt)  -- (\tikztotarget)
    },
  },
}

\newcommand{\phantomarcThroughThreePoints}[4][]{
\coordinate (middle1) at ($(#2)!.5!(#3)$);
\coordinate (middle2) at ($(#3)!.5!(#4)$);
\coordinate (aux1) at ($(middle1)!1!90:(#3)$);
\coordinate (aux2) at ($(middle2)!1!90:(#4)$);
\coordinate (center#2) at ($(intersection of middle1--aux1 and middle2--aux2)$);
\draw[#1,draw=none,name path=kamaan#2] %kamaan#2 is the name of the arc that is being drawn, and center of the arc is center#1 
let \p1=($(#2)-(center#2)$),
      \p2=($(#4)-(center#2)$),
      \n0={veclen(\p1)},       % Radius
      \n1={atan2(\y1,\x1)}, % angles
      \n2={atan2(\y2,\x2)},
      \n3={\n2>\n1?\n2:\n2+360}
    in (#2) arc(\n1:\n3:\n0);
}

\newcommand{\arcThroughThreePoints}[3]{
\phantomarcThroughThreePoints{#1}{#2}{#3}; 
\draw (#1) to[anticlockwise arc centered at=center#1] (#3); %we could give different options to draw for drawing thick arc, colored arc, etc.
}

\keywords{hyperplane arrangement, M\"obius function, characteristic polynomial, Tits algebra, characteristic element, intrinsic volumes}

\begin{document}

\title{Characteristic elements for real hyperplane arrangements}

\author[M.~Aguiar]{Marcelo Aguiar}
\address{Department of Mathematics\\
Cornell University\\
Ithaca, NY 14853}
\email{maguiar@math.cornell.edu}
\urladdr{http://www.math.cornell.edu/~maguiar}

\author[J.~Bastidas]{Jose Bastidas}
\address{Department of Mathematics\\
Cornell University\\
Ithaca, NY 14853}
\email{jdb394@cornell.edu}

\author[S.~Mahajan]{Swapneel Mahajan}
\address{Department of Mathematics\\
Indian Institute of Technology Mumbai\\
Powai, Mumbai 400 076\\ India}
\email{swapneel@math.iitb.ac.in}
\urladdr{http://www.math.iitb.ac.in/~swapneel}

\begin{abstract}
Characteristic elements of the Tits algebra of a real hyperplane arrangement carry  information about the characteristic polynomial. We present this notion and its basic properties, and apply it to derive various results about the characteristic polynomial of an arrangement, from Zaslavsky's formulas to more recent results of Kung and of Klivans and Swartz. We construct several examples of characteristic elements, including one in terms of intrinsic volumes of faces of the arrangement.
\end{abstract}

\maketitle

\section*{Introduction}

We further develop the theory of \emph{characteristic elements} for real hyperplane arrangements started in \cite[Chapter 12]{am17}. These elements of the Tits algebra determine the characteristic polynomial of the arrangement and also determine the characteristic polynomial of the arrangements under each flat. They are defined by requiring that the simple characters of the Tits algebra evaluate on a characteristic element to powers of a specified parameter.

Faces, flats, and the Tits algebra are briefly reviewed in Section \ref{s:tits}.  The fact that the characteristic polynomial of an arrangement, which is defined in terms of flats, carries information about the decomposition of space into faces, originates in work of Zaslavsky  \cite{zaslavsky}, and is at the root of the combinatorial theory of hyperplane arrangements \cite{greenezas,stanley2007}. The Tits algebra provides a natural setting in which this connection can be further gleaned and developed. 
Each arrangement possesses many characteristic elements, and the interest is in constructing particular elements from which specific information about the characteristic polynomial can be extracted. This paper illustrates this fact repeatedly.

We review the notion of characteristic elements in Section \ref{s:char}, extending the definitions and main results of  \cite[Section 12.4]{am17} from linear to affine arrangements. The first applications
are given in Section \ref{s:applications}: we derive the fundamental recursion for the characteristic polynomial from a basic functoriality property of characteristic elements, and we employ multiplicativity of characteristic elements to derive an interesting identity 
% for the characteristic polynomial on a product of two variables
due to Kung \cite{kung04}.
Certain characteristic elements of parameters $1$ and $-1$ are discussed in Section \ref{s:unit-takeuchi}, and employed to derive Zaslavsky's formulas. Section \ref{s:adams} builds
characteristic elements for the simplest Coxeter arrangements in terms of lattice point counting. A main contribution of this paper is the construction of a characteristic element
canonically associated to each arrangement in terms of intrinsic volumes. This is done in Section \ref{s:intrinsic}. As an application, we derive a beautiful result of Klivans and Swartz \cite{KlivansSwartz} which relates the coefficients of the characteristic polynomial to the intrinsic volumes of the chambers.
Characteristic elements can be combined with techniques from the theory of Hopf monoids in species \cite{am10} to study exponential sequences of arrangements, in the sense of Stanley \cite{Stanley96}. This last point will be developed in a longer version of this work. 

\section{The Tits algebra}\label{s:tits}

Let $\arr$ be a (real, affine, finite) hyperplane arrangement: a finite collection of affine hyperplanes in a finite-dimensional real vector space $V$. The hyperplanes in $\arr$ split $V$ into a collection  $\face[\arr]$ of convex polyhedra called \emph{faces}. Given two faces $F$ and $G$, there is a unique face $FG$ containing both $F$ and a small ray from a point in the relative interior of $F$ to a point in the relative interior of $G$.
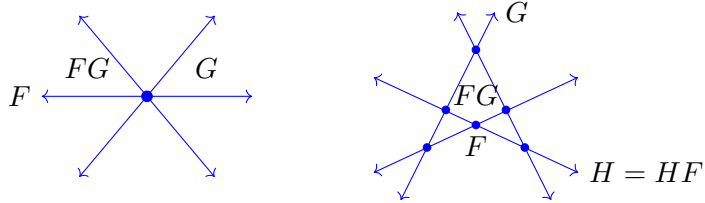
\begin{align*}
&
\begin{gathered}
\begin{tikzpicture}[scale=.7,<->,blue]
\draw (0+180:2) -- (0,0) -- (0:2); 
\draw (50+180:2) -- (0,0) -- (50:2); %Draw a line at 50 degrees and of length 2
\draw (130+180:2) -- (0,0) -- (130:2);
%\draw [draw=none, fill=blue!20!white,fill opacity=.5] (0,0) -- (130:1.8) -- (180:1.8) -- cycle;
 \draw [fill=blue] (0,0) circle (.1);
\node [left,black] at (0+180:2) {\small $F$};
\node [black] at (25:1.25) {\small $G$}; %placing the label exactly in the middle of the two rays
\node [black] at (155:1.25) {\small $FG$};
\end{tikzpicture} 
\end{gathered}
& &
\begin{gathered}
\begin{tikzpicture}[scale=.5,<->,blue]
\draw (25+180:3) -- (0,0) -- (25:3); 
\draw (155+180:3) -- (0,0) -- (155:3); 
\draw   (-2,-2) -- (0.5,3); 
\draw   (2,-2) -- (-0.5,3); 
\draw [fill=blue] (0,0) circle (.1);
\draw [fill=blue] (0,2) circle (.1);
\draw [fill=blue] (.8,0.4) circle (.1);
\draw [fill=blue] (-.8,0.4) circle (.1);
\draw [fill=blue] (1.3,-0.6) circle (.1);
\draw [fill=blue] (-1.3,-0.6) circle (.1);
\node [below,black] at (0,0) {\small $F$};
\node [right,black] at (0.5,3) {\small $G$};
\node [above,right,black] at (155+180:3) {\small $H=HF$};
\node [black] at (0,0.8) {\small $FG$};
\end{tikzpicture} 
\end{gathered}
\end{align*}
Endowed with this operation, the set $\face[\arr]$ is a semigroup.
Its linearization over a field $\Kb$ is the \emph{Tits algebra} $\Kb\face[\arr]$.
See \cite[Chapters 1 and 9]{am17} for more details.
For linear arrangements, such as the one on the left above, $\face[\arr]$ is a monoid and the central face is the unit. For affine arrangements, such as the one on the right, this semigroup is nonunital.
An interesting fact that we review below (Theorem \ref{t:unit-affine}) is that the Tits algebra is always unital. We let $\B{H}_F$ denote the basis element of $\Kb\face[\arr]$ associated to the face $F$ of $\arr$.

An arbitrary intersection of hyperplanes is a \emph{flat} of the arrangement. The set $\flat[\arr]$ of flats is a join-semilattice. The support of a face $F$ is the smallest flat $\supp(F)$ containing it. 
We view $\flat[\arr]$ as a commutative semigroup with the join operation for the product and then the support map
\[
\supp:\face[\arr]\to\flat[\arr]
\]
is a morphism of semigroups.

The algebra $\Kb \flat[\arr]$ is also unital and is the maximal semisimple quotient of $\Kb\face[\arr]$ via the support map. We let $\B{H}_\rX$ denote the basis element of $\Kb \flat[\arr]$ associated to the flat $\rX$ of $\arr$. It follows that the simple modules over the Tits algebra are one-dimensional and are indexed by flats. The character $\chi_\rX$ associated to the flat $\rX$ evaluated on an element 
\begin{equation}\label{eq:elt-tits}
w = \sum_F w^F \B{H}_F
\end{equation}
of $\Kb\face[\arr]$ yields 
\begin{equation}\label{eq:elt-char}
\chi_\rX(w) = \sum_{F:\,\supp(F)\leq\rX} w^F.
\end{equation}
This rests on the fact that the unique complete system of orthogonal idempotents for $\Kb \flat[\arr]$ consists of elements $\B{Q}_\rX$ uniquely determined by
\[
\B{H}_\rY = \sum_{\rX:\,\rY\leq\rX}\B{Q}_\rX.
\]
This is a result of Solomon \cite[Theorem 1]{solomon67}.

The sets $\face[\arr]$ and $\flat[\arr]$ are partially ordered by inclusion. Both posets are graded and of the same rank $\rank(\arr)$, the rank of the arrangement. The maximal elements of $\face[\arr]$ are the \emph{chambers}. The ambient space is the top element of $\flat[\arr]$, we denote it by $\maxflat$. 

The \emph{characteristic polynomial} of $\arr$ is
\begin{equation}\label{eq:char-poly}
\cpoly(\arr,t) := \sum_{\rY} \,\mu(\rY,\maxflat) \,t^{\rank(\rY)}. 
\end{equation}
The sum is over all flats. It is a monic polynomial of degree $\rank(\arr)$.

The \emph{arrangement under a flat $\rX$} is
\[
\arr^\rX = \{\Hr\cap\rX \mid \Hr\in\arr,\, \rX\not\subseteq\Hr,\, \Hr\cap\rX\neq\emptyset\}.
\]
The flats of $\arr^\rX$ are the flats of $\arr$ that are contained in $\rX$. Hence,
\begin{equation}\label{eq:char-poly-under}
\cpoly(\arr^\rX,t) := \sum_{\rY:\,\rY\leq\rX} \,\mu(\rY,\rX) \,t^{\rank(\rY)}. 
\end{equation}

\section{Characteristic elements}\label{s:char}

The definitions and results in this section extend those of \cite[Section 12.4]{am17} to the setting of affine arrangements. 

\subsection{Definition and basic properties}\label{ss:def-basic}

Let $t$ be a fixed scalar.
An element $w$ of the Tits algebra is \emph{characteristic of parameter $t$} if for each flat $\rX$
\begin{equation}\label{eq:char-def}
\chi_\rX(w) = t^{\rank(\rX)},
\end{equation}
with $\chi_\rX(w)$ as in \eqref{eq:elt-char}. 

Two characteristic elements of the same parameter take the same value on all simple modules, and hence differ by a nilpotent element (an element of the Jacobson radical). The set of characteristic elements of a given parameter is an affine subspace of the Tits algebra of dimension equal to the number of faces minus the number of flats. 

One-dimensional characters are multiplicative. We deduce the following result. 
%\cite[Lemma 12.54]{am17}

\begin{lemma}\label{l:char-mult}
If $u$ is a characteristic element of parameter $s$ and
$v$ is a characteristic element of parameter $t$, then
$uv$ is a characteristic element of parameter $st$. 
\end{lemma}

\subsection{Relation to the characteristic polynomial}\label{ss:char-cpoly}

The right-hand sides of \eqref{eq:char-poly-under} and \eqref{eq:char-def} are related by M\"obius inversion, which implies the following result.
% \cite[Lemma 12.55]{am17}

\begin{lemma}\label{l:char-cpoly}
An element $w$ of the Tits algebra is characteristic of parameter $t$ iff for every flat $\rX$,
\begin{equation}\label{eq:char-cpoly-under}
\sum_{F:\,\supp(F)=\rX} w^F = \cpoly(\arr^\rX,t).
\end{equation}
\end{lemma}
In particular, since the chambers are the faces of top support:
\begin{equation}\label{eq:char-cpoly}
\sum_{C} w^C = \cpoly(\arr,t),
\end{equation}
with the sum over all chambers $C$ of $\arr$.

\subsection{Functoriality}\label{ss:functor}

Let $\arr'$ be a subarrangement of $\arr$: $\arr'$ consists of some of the hyperplanes in $\arr$.
There is a morphism of semigroups 
\begin{equation}\label{eq:tits-to-subarr}
\face[\arr]\to \face[\arr']  
\end{equation}
which sends a face $F$ of $\arr$ to the unique face of $\arr'$ whose interior contains the interior of $F$. This in turn induces a morphism from the Tits algebra of $\arr$ to that of $\arr'$: if $w$ is as in \eqref{eq:elt-tits}, then
\[
f(w)=\sum_{G\in\face[\arr']} f(w)^G\,\B{H}_G, \text{ \ where \ } f(w)^G = \sum_{F:\,f(F)=G} w^F.
\]

%\cite[Lemma 12.63]{am17}  
\begin{lemma}\label{l:char-subarr}
Assume that $\arr$ and $\arr'$ have the same rank. Then $f$ sends characteristic elements for $\arr$ to characteristic elements for $\arr'$, of the same parameter.
\end{lemma}

\section{First applications}\label{s:applications}

\subsection{The fundamental recursion for the characteristic polynomial}\label{ss:removing}

The characteristic polynomial of $\arr$ may be calculated recursively by removing one hyperplane at a time. 
As a first application, we derive a proof of this formula.

Let $\Hr$ be a hyperplane in $\arr$ and set $\arr'=\arr\setminus\{\Hr\}$. Assume that $\rank(\arr')=\rank(\arr)$.

Pick any characteristic element $w$ of parameter $t$. Applying \eqref{eq:char-cpoly-under} to calculate the characteristic polynomial of $\arr^\Hr$, and \eqref{eq:char-cpoly} to calculate that of $\arr$, we obtain
\[
\cpoly(\arr,t) + \cpoly(\arr^\Hr,t) = \sum w^C + \sum w^F.
\]
The first sum is over all chambers $C$ of $\arr$ and the second over all faces $F$ of $\arr$ with $\supp(F)=\Hr$. By Lemma \ref{l:char-subarr}, we may further employ \eqref{eq:char-cpoly} to calculate $\cpoly(\arr',t)$ in terms of coefficients of $f(w)$.  We obtain
\[
\cpoly(\arr',t) = \sum f(w)^D = \sum w^G
\]
the first sum being over all chambers $D$ of $\arr'$ and the second over all faces $G$ of $\arr$ with $f(G)=D$ for some such $D$.  These faces $G$ are either chambers of $\arr$ or faces with support $\Hr$. Comparing the above expressions, we conclude that
\begin{equation}\label{eq:removing-recursion}
\cpoly(\arr,t) = \cpoly(\arr\setminus\{\Hr\},t) - \cpoly(\arr^\Hr,t).
\end{equation}

We gave this derivation of the fundamental recursion in \cite[Proposition 12.66]{am17} (for linear arrangements). The proof in \cite[Lemma 2.2]{stanley2007}, \cite[Theorem 2.56]{orlikterao} is quite different.  

\subsection{The characteristic polynomial on a product. An identity of Kung}\label{ss:product}

For the second application we employ Lemma \ref{l:char-mult}. Pick characteristic elements $u$ and $v$ of parameters $s$ and $t$, respectively. Applying \eqref{eq:char-cpoly} to the characteristic element $uv$, we obtain
\[
\cpoly(\arr,st) = \sum_C (uv)^C = \sum_{C,F,G:\, FG=C} u^F v^G.
\]
The first sum is over chambers $C$ and the second over faces $F$ and $G$ which multiply to a chamber. This happens precisely when 
$\supp(F)\join\supp(G) = \supp(FG)= \maxflat$,
since $\supp$ is a morphism of semigroups. So the previous sum equals
\[
\sum_{\rX,\rY:\,\rX\join\rY=\maxflat} \sum_{\substack{F:\,\supp(F)=\rX\\G:\,\supp(G)=\rY}} u^F v^G.
\]
Combining the preceding with \eqref{eq:char-poly-under} we obtain
\[
\cpoly(\arr,st) = \sum_{\rX,\rY:\,\rX\join\rY=\maxflat} \cpoly(\arr^\rX,s) \cpoly(\arr^\rY,t).
\]
Finally, an application of \cite[Lemma 1.86]{am17} yields
\begin{equation}\label{eq:char-kung}
\cpoly(\arr,st) = \sum_{\rX} t^{\rank(\rX)} \cpoly(\arr^\rX,s) \cpoly(\arr_\rX,t),
\end{equation}
where $\arr_\rX$ denotes the arrangement \emph{over} the flat $\rX$ \cite[Section 1.7.2]{am17}. 
This identity is due to Kung \cite[Theorem 4]{kung04}.
Kung discusses a couple of proofs, all quite different from the one above. Kung's result is for matroids, which covers the case of linear arrangements. The identity above holds for affine arrangements.

\section{Characteristic elements of parameters $\pm 1$ }\label{s:unit-takeuchi}

\subsection{The unit element}

An affine arrangement $\arr$ is \emph{essential} if the minimal flats are points. In general, the minimal flats are pairwise parallel subspaces of a common dimension. Intersecting with an orthogonal subspace makes $\arr$ essential. Faces of $\arr$ are in correspondence with faces of the essentialization. The same applies to flats.~A face of $\arr$ is \emph{essentially bounded} if~the corresponding face of the essentialization is bounded. The following is \cite[Theorem~14.23]{am17}.

\begin{theorem}\label{t:unit-affine}
The Tits algebra of an affine arrangement $\arr$ possesses a unit element. The unit is 
\begin{equation}\label{eq:unit-affine}
\upsilon = \sum_F (-1)^{\rank(F)} \B{H}_F,
\end{equation}
with $F$ running over the set of essentially bounded faces of $\arr$.
\end{theorem}

When $\arr$ is linear, the only essentially bounded face is the central face $\minface$, and $\upsilon=\B{H}_\minface$.

The unit element acts as the identity on any module, and hence the one-dimensional characters evaluate to $1$ on it. This implies the following result. 

\begin{proposition}\label{p:unit-char}
The unit element $\upsilon$ is characteristic of parameter $1$.
\end{proposition}

Here is an alternative proof of the proposition. According to \eqref{eq:elt-char}, the character value
$\chi_\rX(\upsilon) = \sum (-1)^{\rank(F)}$ is the Euler characteristic of the complex consisting of the essentially bounded part of $\rX$, as illustrated below. The latter is contractible \cite[Theorem 4.5.7]{blswz},
% also \cite[Chapter 1, Exercise 7]{stanley07}
and hence the character value is $1$. 
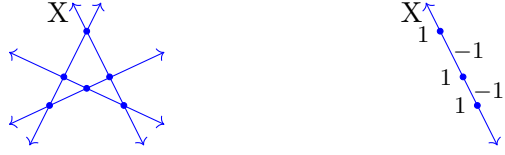
\begin{align*}
&
\begin{gathered}
\begin{tikzpicture}[scale=.38,<->,blue]
\draw (25+180:3) -- (0,0) -- (25:3); 
\draw (155+180:3) -- (0,0) -- (155:3); 
\draw   (-2,-2) -- (0.5,3); 
\draw   (2,-2) -- (-0.5,3); 
\draw [fill=blue] (0,0) circle (.1);
\draw [fill=blue] (0,2) circle (.1);
\draw [fill=blue] (.8,0.4) circle (.1);
\draw [fill=blue] (-.8,0.4) circle (.1);
\draw [fill=blue] (1.3,-0.6) circle (.1);
\draw [fill=blue] (-1.3,-0.6) circle (.1);
\node [left,below,black] at (-1,3.4) {\small $\rX$};
\end{tikzpicture} 
\end{gathered}
& &
\begin{gathered}
\begin{tikzpicture}[scale=.38,<->,blue]
\draw   (2,-2) -- (-0.5,3); 
\draw [fill=blue] (0,2) circle (.1);
\draw [fill=blue] (.8,0.4) circle (.1);
\draw [fill=blue] (1.3,-0.6) circle (.1);
\node [left,black] at (0,1.9) {\scriptsize $1$};
\node [left,black] at (.8,0.4) {\scriptsize $1$};
\node [left,black] at (1.3,-0.6) {\scriptsize $1$};
\node [black] at (1,1.3) {\scriptsize $-1$};
\node [black] at (1.7,-0.1) {\scriptsize $-1$};
\node [left,below,black] at (-1,3.4) {\small $\rX$};
\end{tikzpicture} 
\end{gathered}
\end{align*}
The proof of Theorem \ref{t:unit-affine} is also topological, but more involved.

\subsection{The Takeuchi element}

The \emph{Takeuchi element} is
\begin{equation}\label{eq:takeuchi}
\tau = \sum_F \,(-1)^{\rank(F)}\,\B{H}_F, 
\end{equation}
with the sum over all faces $F$ of $\arr$.

The following extends \cite[Corollary 12.57]{am17} to the setting of affine arrangements.

\begin{proposition}\label{p:tak-char}
The Takeuchi element $\tau$ is characteristic of parameter $-1$.
\end{proposition}

This time the proof can be brought down to the calculation of the Euler characteristic of a relative pair of cell complexes $(B,\partial B)$, where $B$ is the complex obtained by dissecting a large ball (containing the bounded faces) with the hyperplanes in $\arr$.

\subsection{Application: Zaslavsky's formulas}\label{ss:Zaslavsky}

All chambers of $\arr$ are of rank $\rank(\arr)$. Applying \eqref{eq:char-cpoly} to the unit element $\upsilon$ we obtain that
\[
(-1)^{\rank(\arr)} \cpoly(\arr,1) = (-1)^{\rank(\arr)} \sum_\rY \mu(\rY,\maxflat)
\]
equals the number of essentially bounded chambers in $\arr$. Employing the Takeuchi element $\tau$ instead, we obtain that
\[
(-1)^{\rank(\arr)} \cpoly(\arr,-1) = \sum_\rY (-1)^{\rank(\arr)-\rank(\rY)} \mu(\rY,\maxflat)
\]
equals the total number of chambers in $\arr$. These are Zaslavsky's formulas
\cite[Theorem A, Theorem C, Corollary 2.2]{zaslavsky}, \cite[Proposition 8.1]{lasvergnas75}.

\begin{remark}
The above proof does not differ substantially from Zaslavsky's. The core topological argument has been shifted to prove the facts that $\upsilon$ and $\tau$ are characteristic.
\end{remark}

\section{The Adams elements}\label{s:adams}

\subsection{Braid arrangement}\label{ss:adamsA} 

The \emph{braid arrangement} $\arr_n$ in $\Rb^n$ consists of the diagonal hyperplanes
$
x_i=x_j
$
for $1\leq i<j\leq n$. It is a linear arrangement of rank $n-1$. The central face is the line $x_1=\cdots=x_n$. Intersecting with a sphere around the origin in the hyperplane $x_1+\cdots+x_n=0$ we obtain the \emph{Coxeter complex of type $A_{n-1}$}. The pictures below show the cases $n=3$ and $4$.
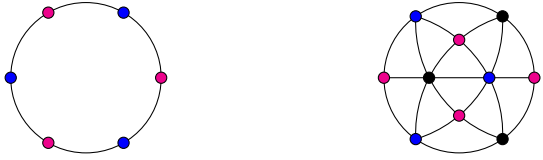
\begin{align*}
&
\begin{gathered}
\begin{tikzpicture}[scale=.4]
\newdimen\R
\R=2.5cm %radius of the outer circle
\coordinate (madhya) at (0,0);
\draw (madhya) circle (\R); %drawing the outer circle
\coordinate (A1) at (60:\R);
\coordinate (B1) at (180+60:\R);
\coordinate (A2) at (180-60:\R);
\coordinate (B2) at (-60:\R);
\coordinate (A3) at (0:\R);
\coordinate (B3) at (180:\R);
\node [circle, inner sep=1.5pt,fill=blue,draw] at (A1) {};
\node [circle, inner sep=1.5pt,fill=magenta,draw] at (B1) {};
\node [circle, inner sep=1.5pt,fill=magenta,draw] at (A2) {};
\node [circle, inner sep=1.5pt,fill=blue,draw] at (B2) {};
\node [circle, inner sep=1.5pt,fill=magenta,draw] at (A3) {};
\node [circle, inner sep=1.5pt,fill=blue,draw] at (B3) {};
\end{tikzpicture}
\end{gathered}
& &
\begin{gathered}
\begin{tikzpicture}[scale=.4]
\newdimen\R
\R=2.5cm %radius of the outer circle
\coordinate (madhya) at (0,0);
\draw (madhya) circle (\R); %drawing the outer circle
\coordinate (P) at (-1,0); %the first vertex
\coordinate (Q) at (1,0); %the second vertex
\coordinate (A1) at (54.735:\R);
\coordinate (B1) at (180+54.735:\R);
\arcThroughThreePoints{A1}{P}{B1}; 
\arcThroughThreePoints{B1}{Q}{A1}; 
\coordinate (A2) at (180-54.735:\R);
\coordinate (B2) at (-54.735:\R);
\arcThroughThreePoints{A2}{P}{B2};
\arcThroughThreePoints{B2}{Q}{A2};
\coordinate (A3) at (0:\R);
\coordinate (B3) at (180:\R);
\draw (A3) -- (B3);
\path[name intersections={of=kamaanA1 and kamaanB2, by=a1b2}]; 
\path[name intersections={of=kamaanB1 and kamaanA2, by=b1a2}]; 
\node [circle, inner sep=1.5pt,fill=magenta,draw] at (a1b2) {};
\node [circle, inner sep=1.5pt,fill=magenta,draw] at (b1a2) {};
\node [circle, inner sep=1.5pt,fill=black,draw] at (P) {};
\node [circle, inner sep=1.5pt,fill=blue,draw] at (Q) {};
\node [circle, inner sep=1.5pt,fill=black,draw] at (A1) {};
\node [circle, inner sep=1.5pt,fill=blue,draw] at (B1) {};
\node [circle, inner sep=1.5pt,fill=blue,draw] at (A2) {};
\node [circle, inner sep=1.5pt,fill=black,draw] at (B2) {};
\node [circle, inner sep=1.5pt,fill=magenta,draw] at (A3) {};
\node [circle, inner sep=1.5pt,fill=magenta,draw] at (B3) {};
\end{tikzpicture}
\end{gathered}
\end{align*}

Given a face $F$ of $\arr_n$, let $\deg(F)$ denote its dimension as a subset of $\Rb^n$ (so two more than the dimension as a subset of the above sphere). The rank of $F$ in the face poset is then $\deg(F)-1$. Analogous considerations apply to flats.

The \emph{Adams element of type $A_{n-1}$} (and parameter $t$) is defined by
\begin{equation}\label{eq:adamsA}
\alpha_t = \sum_F  \binom{t}{\deg(F)} \B{H}_F,
\end{equation}
with the sum over the faces $F$ of $\arr_n$. For each integer $k$, the binomial coefficient $\binom{k}{\deg(F)}$ counts the number of points in the relative interior of $F$ with coordinates from $[k]=\{1,\dots,k\}$. On the other hand, given a flat $\rX$, the number $k^{\deg(\rX)}$ is the number of points in $\rX\cap [k]^n$. Since $\rX$ splits as the disjoint union of the relatively open faces $F$ with $\supp(F)\leq\rX$, we have that
\[
\sum_{F:\,\supp(F)\leq\rX} \binom{k}{\deg(F)} = k^{\deg(\rX)}.
\]
We have shown the following, for which a different proof is given in \cite[Lemma 12.78]{am17}.

\begin{proposition}\label{p:adamsA}
For any nonzero scalar $t$, the element $\frac{1}{t}\alpha_t$ is characteristic of parameter $t$.
\end{proposition}

There are $n!$ chambers in $\arr_n$. As a small application of \eqref{eq:char-cpoly}, we obtain the well-known expression for the characteristic polynomial of the braid arrangement.

\begin{equation}\label{eq:braid-cpoly}
\cpoly(\arr_n,t)= \frac{1}{t} \sum_{C} \binom{t}{n} = (t-1)(t-2)\cdots(t-(n-1)).
\end{equation}

It follows from Lemma \ref{l:char-mult} that $\frac{1}{st}\alpha_s\alpha_t$ is characteristic of parameter $st$. In fact, it can be shown that $\alpha_s\alpha_t=\alpha_{st}$, see for instance \cite[Lemma 12.80]{am17}.

\subsection{Signed braid arrangement}\label{ss:adamsB}

The \emph{signed braid arrangement} $\arr^{\pm}_n$  in $\Rb^n$ consists of the hyperplanes 
$x_i=x_j$, $x_i=-x_j$ for $1 \leq i<j \leq n$ and $x_k=0$ for $1 \leq k \leq n$. 
It has rank $n$. The \emph{Coxeter complex of type $B_n$}, obtained by intersecting with the sphere in $\Rb^n$, is shown below for $n=2$ and $3$. 

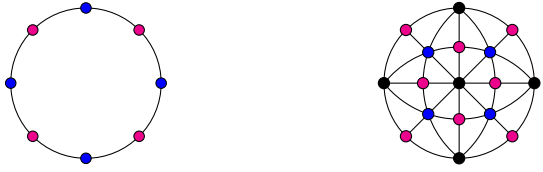
\begin{align*}
\begin{gathered} 
\begin{tikzpicture}[scale=1]
\draw circle (1);
\node [circle, inner sep=1.4pt,fill=blue,draw] at (0:1) {};
\node [circle, inner sep=1.4pt,fill=magenta,draw] at (45:1) {};
\node [circle, inner sep=1.4pt,fill=blue,draw] at (90:1) {};
\node [circle, inner sep=1.4pt,fill=magenta,draw] at (135:1) {};
\node [circle, inner sep=1.4pt,fill=blue,draw] at (180:1) {};
\node [circle, inner sep=1.4pt,fill=magenta,draw] at (225:1) {};
\node [circle, inner sep=1.4pt,fill=blue,draw] at (270:1) {};
\node [circle, inner sep=1.4pt,fill=magenta,draw] at (315:1) {};
\end{tikzpicture} 
\end{gathered} 
& &
\begin{gathered}
\begin{tikzpicture}[scale=0.4]
\newdimen\R
\R=2.5cm %radius of the outer circle
\coordinate (madhya) at (0,0);
\draw (madhya) circle (\R); %drawing the outer circle
\coordinate (P) at (-1.2,0); 
\coordinate (Q) at (1.2,0); 
\coordinate (R) at (0,1.2); 
\coordinate (S) at (0,-1.2); 
\coordinate (A1) at (0:\R);
\coordinate (B1) at (180:\R);
\draw (A1) -- (B1);
\arcThroughThreePoints{A1}{R}{B1}; 
\arcThroughThreePoints{B1}{S}{A1}; 
\coordinate (A2) at (45:\R);
\coordinate (B2) at (225:\R);
\draw (A2) -- (B2);
\coordinate (A3) at (90:\R);
\coordinate (B3) at (270:\R);
\draw (A3) -- (B3);
\arcThroughThreePoints{A3}{P}{B3}; 
\arcThroughThreePoints{B3}{Q}{A3}; 
\coordinate (A4) at (135:\R);
\coordinate (B4) at (315:\R);
\draw (A4) -- (B4);
\path[name intersections={of=kamaanA1 and kamaanB3, by=a1b3}]; 
\path[name intersections={of=kamaanB1 and kamaanA3, by=b1a3}]; 
\path[name intersections={of=kamaanA1 and kamaanA3, by=a1a3}]; 
\path[name intersections={of=kamaanB1 and kamaanB3, by=b1b3}]; 
\node [circle, inner sep=1.5pt,fill=black,draw] at (madhya) {};
\node [circle, inner sep=1.5pt,fill=blue,draw] at (a1b3) {};
\node [circle, inner sep=1.5pt,fill=blue,draw] at (b1a3) {};
\node [circle, inner sep=1.5pt,fill=blue,draw] at (a1a3) {};
\node [circle, inner sep=1.5pt,fill=blue,draw] at (b1b3) {};
\node [circle, inner sep=1.5pt,fill=magenta,draw] at (P) {};
\node [circle, inner sep=1.5pt,fill=magenta,draw] at (Q) {};
\node [circle, inner sep=1.5pt,fill=magenta,draw] at (R) {};
\node [circle, inner sep=1.5pt,fill=magenta,draw] at (S) {};
\node [circle, inner sep=1.5pt,fill=black,draw] at (A1) {};
\node [circle, inner sep=1.5pt,fill=black,draw] at (B1) {};
\node [circle, inner sep=1.5pt,fill=magenta,draw] at (A2) {};
\node [circle, inner sep=1.5pt,fill=magenta,draw] at (B2) {};
\node [circle, inner sep=1.5pt,fill=black,draw] at (A3) {};
\node [circle, inner sep=1.5pt,fill=black,draw] at (B3) {};
\node [circle, inner sep=1.5pt,fill=magenta,draw] at (A4) {};
\node [circle, inner sep=1.5pt,fill=magenta,draw] at (B4) {};
\end{tikzpicture}
\end{gathered} 
\end{align*}

The \emph{Adams element of type $B_n$} (and odd parameter $2t+1$) is defined by
\begin{equation}\label{eq:adamsB}
\alpha^{\pm}_{2t+1} = \sum_F \,\binom{t}{\rank(F)} \,\B{H}_F.
\end{equation}
Proceeding as in Section \ref{ss:adamsA}, but now counting integer points in $[-k,k]^n\cap\rX$ for each flat $\rX$ according to the face in which they lie, one arrives at the following fact, for which a different proof is given in \cite[Proposition 12.89]{am17}.

\begin{proposition}\label{p:adamsB}
For any scalar $t$, the element $\alpha^{\pm}_{2t+1}$ is characteristic of parameter $2t+1$.
\end{proposition}

There are $(2n)!!$ chambers in $\arr^{\pm}_n$. Employing \eqref{eq:char-cpoly}, one obtains that
\[
\cpoly(\arr^{\pm}_n,2t+1)=  (2n)!! \binom{t}{n},
\]
which is equivalent to the familiar expression for the characteristic polynomial of the signed braid arrangement:
\begin{equation}\label{eq:signedbraid-cpoly}
\cpoly(\arr^{\pm}_n,t)=   (t-1)(t-3)\cdots(t-(2n-1)).
\end{equation}

Related elements $\alpha^{\pm}_{2t}$ are discussed in \cite[Section 12.6.3]{am17}. These are not characteristic.

\subsection{Coordinate arrangement}\label{ss:adams-coord}

The coordinate arrangement $\carr_n$ in $\Rb^n$ consists of the hyperplanes
$
x_i=0
$
for $1\leq i\leq n$. The associated subdivision of the sphere is the \emph{Coxeter complex of type $A_1^n$}. The \emph{first orthant} is $\bigcap_i\{x_i\geq 0\}$. For each face $F$ of $\carr_n$, let 
\[
\gamma_t^F = \begin{cases}
(t-1)^{\rank(F)} & \text{ if $F$ lies in the first orthant,} \\
 0        & \text{ if not.}
\end{cases}
\]
An argument similar to those in Sections \ref{ss:adamsA} and \ref{ss:adamsB}, but now counting integer points in $[0,k-1]^n\cap\rX$, shows that the element 
\[
\gamma_t = \sum_F \gamma_t^F \B{H}_F
\]
is characteristic of parameter $t$. In this case, only one chamber appears with nonzero coefficient in $\gamma_t$ (the first orthant). We obtain 
\[
\cpoly(\carr_n,t) = (t-1)^n.
\]

\begin{remark}
The strategy employed in this section to build characteristic elements draws on ideas of Beck and Zaslavsky in \cite{bz06}, and in fact may be further developed to study the polynomials introduced in that work.
\end{remark}

\section{Intrinsic elements}\label{s:intrinsic}

We employ the notion of \emph{intrinsic volumes} of a convex cone $C$ in $\Rb^n$ \cite[Section 2.2]{amelunxen15}. For each $k=0,\dots,n$, let $v_k(C)$ be the proportion of volume of space occupied by points that map to a $k$-dimensional face of $C$ under the \emph{nearest point} projection. 
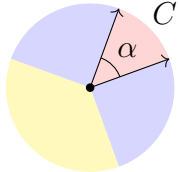
\begin{align*}
 \begin{gathered}
 \begin{tikzpicture}[scale=.7]
 \newdimen\R
 \R=2cm %radius of the outer circle
 \coordinate (center) at (0,0); 
 \draw [draw=none, fill=red!20!white,fill opacity=.8] (center) -- (20:.8*\R) arc (20:70:.8*\R) -- cycle;
 \draw [draw=none, fill=blue!20!white,fill opacity=.8] (center) -- (70:.8*\R) arc (70:70+90:.8*\R) -- (center);
 \draw [draw=none, fill=blue!20!white,fill opacity=.8] (center) -- (20:.8*\R) arc (20:20-90:.8*\R) -- (center);
 \draw [draw=none, fill=yellow!40!white,fill opacity=.8] (center) -- (70+90:.8*\R) arc (70+90:290:.8*\R)  -- (center);
  \draw[<->]  (20:.8*\R)--(0,0) -- (70:.8*\R); 
  \node [circle, inner sep=1pt,fill=black,draw] at (center) {};
\draw (20:.3*\R) to[anticlockwise arc centered at=center] (70:.3*\R);
\node at (45:.5*\R) {$\alpha$}; 
%\node at (70+45:.5*\R) {\small $1/4$}; \node at (20-45:.5*\R) {\small $1/4$};
\node at (45:\R) {$C$};
 \end{tikzpicture} 
 \end{gathered}
 & & 
 \begin{gathered} 
 v_2(C)=\alpha/2\pi\\
 v_1(C)=1/2\\
 v_0(C)=1/2-\alpha/2\pi
 \end{gathered}
 \end{align*}

We record the following properties. Let $C$ be a cone. If $k>\dim(C)$, or if $k$ is smaller than the dimension of the minimal face of $C$, then $v_k(C)=0$.
Also,
 \begin{equation}\label{eq:intrinsic-distribution}
  \sum_{k=0}^n v_k(C) = 1.
\end{equation}
Most importantly, each intrinsic volume $v_k$ is a \emph{valuation} on convex cones.
The \emph{Gauss-Bonnet formula} states that if $C$ is not a subspace, then
 \begin{equation}\label{eq:gauss-bonnet}
 \sum_{k=0}^{n} (-1)^k v_k(C) =  0.
 \end{equation} 
If $C$ is a subspace, then
 \begin{equation}\label{eq:intrinsic-subspace}
 v_k(C) = \begin{cases}
1 & \text{ if } \dim(C)=k,\\
 0 & \text{ if not.}
\end{cases}
 \end{equation} 
 
We extend this notion to convex polyhedra $P$, with the same definition. It then turns out that $v_k(P)$ depends only on the \emph{recession cone} of $P$. Each $v_k$ is a valuation on convex polyhedra.
  
Let $\arr$ be an arrangement in $\Rb^n$ and $d$ the dimension of the minimal faces of $\arr$.
Each face of $\arr$ is a polyhedron.
We define the \emph{intrinsic element} of parameter $t$ for $\arr$ by
\begin{equation}\label{eq:intrinsic-elt}
\nu_t = \sum_F (-1)^{\dim(F)} \Bigl(\sum_{k=d} ^{\dim(F)} (-1)^{k} v_k(F) t^{k-d}\Bigr) \B{H}_F.
\end{equation}
For any face $F$, $\rank(F)=\dim(F)-d$.

\begin{theorem}\label{t:intrinsic}
The element $\nu_t$ is characteristic of parameter $t$.
\end{theorem}
The proof relies on the valuation property of intrinsic volumes,
plus \eqref{eq:intrinsic-subspace}.

As an immediate consequence of the theorem, we deduce the following result of Klivans and Swartz \cite[Theorem 5]{KlivansSwartz}.

\begin{corollary}\label{c:intrinsic-ks}
The coefficient of $t^j$ in the characteristic polynomial of $\arr$ is
\begin{equation}\label{eq:intrinsic-cpoly}
(-1)^{\rank(\arr)-j} \sum_C v_{j+d}(C),
\end{equation}
with the sum over all chambers of $\arr$.
\end{corollary}

We turn to general properties of the intrinsic elements. 

\begin{proposition}\label{p:intrinsic-mult}
For any parameters $s$ and $t$, $\nu_{s}\nu_{t}=\nu_{st}$.
\end{proposition}
The proof relies on results of McMullen \cite[Theorems 2 and 3]{mcmullen75}.

The recession cone of a face $F$ is a subspace if and only if $F$ is essentially bounded.
Together with the Gauss-Bonnet formula \eqref{eq:gauss-bonnet}, this implies that the intrinsic element of parameter $1$ is  precisely the unit element of the Tits algebra: $\nu_1=\upsilon$. Formula \eqref{eq:intrinsic-distribution} implies that the intrinsic element of parameter
$-1$ and the Takeuchi element coincide: $\nu_{-1}=\tau$.

The following table shows the intrinsic volumes for the braid arrangement of rank $3$ in its canonical realization. The ambient space is the hyperplane $x_1+\cdots+x_4=0$ in $\Rb^4$.
All faces of the same \emph{type} are congruent and hence have the same volumes.
The edges of types $(2,1,1)$ and $(1,1,2)$ have size $x$ and the edges of type $(1,2,1)$ have size $y$, where
\[
x = \frac{1}{2\pi} \arccos(\frac{\sqrt{3}}{3}) \qand y = \frac{1}{2\pi} \arccos(\frac{1}{3}).
\]
The entries above the main diagonal are $0$.
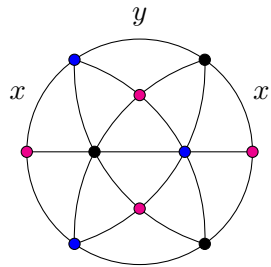
\begin{align*}
&
\begin{gathered}
\begin{tabular}[htp]{|r|c|c|c|c|}
\hline
face  & $v_0$ & $v_1$ & $v_2$ & $v_3$ \\
\hline
center & $1$ &  &  &  \\
vertices & $1/2$ & $1/2$ &  &  \\
short edges  & $1/2-x$  & $1/2$ & $x$ &  \\
long edges  & $1/2-y$ & $1/2$ & $y$ & \\
triangles & $1/4$ & $11/24$ & $1/4$ & $1/24$  \\
\hline
\end{tabular}
\end{gathered}
& &
\begin{gathered}
\begin{tikzpicture}[scale=.6]
\newdimen\R
\R=2.5cm %radius of the outer circle
\coordinate (madhya) at (0,0);
\draw (madhya) circle (\R); %drawing the outer circle
\coordinate (P) at (-1,0); %the first vertex
\coordinate (Q) at (1,0); %the second vertex
\coordinate (A1) at (54.735:\R);
\coordinate (B1) at (180+54.735:\R);
\arcThroughThreePoints{A1}{P}{B1}; 
\arcThroughThreePoints{B1}{Q}{A1}; 
\coordinate (A2) at (180-54.735:\R);
\coordinate (B2) at (-54.735:\R);
\arcThroughThreePoints{A2}{P}{B2};
\arcThroughThreePoints{B2}{Q}{A2};
\coordinate (A3) at (0:\R);
\coordinate (B3) at (180:\R);
\draw (A3) -- (B3);
\path[name intersections={of=kamaanA1 and kamaanB2, by=a1b2}]; 
\path[name intersections={of=kamaanB1 and kamaanA2, by=b1a2}]; 
\node [circle, inner sep=1.5pt,fill=magenta,draw] at (a1b2) {};
\node [circle, inner sep=1.5pt,fill=magenta,draw] at (b1a2) {};
\node [circle, inner sep=1.5pt,fill=black,draw] at (P) {};
\node [circle, inner sep=1.5pt,fill=blue,draw] at (Q) {};
\node [circle, inner sep=1.5pt,fill=black,draw] at (A1) {};
\node [circle, inner sep=1.5pt,fill=blue,draw] at (B1) {};
\node [circle, inner sep=1.5pt,fill=blue,draw] at (A2) {};
\node [circle, inner sep=1.5pt,fill=black,draw] at (B2) {};
\node [circle, inner sep=1.5pt,fill=magenta,draw] at (A3) {};
\node [circle, inner sep=1.5pt,fill=magenta,draw] at (B3) {};
\node at (0,3)  {\small $y$};
\node at (2.7,1.3) {\small $x$};
\node at (-2.7,1.3) {\small $x$};
\end{tikzpicture}
\end{gathered}
\end{align*}

Each circle is composed of four edges of size $x$ and two edges of size $y$, so $4x+2y=1$. The arrangement under a circle is combinatorially isomorphic to the braid arrangement of rank $2$.
Employing \eqref{eq:intrinsic-cpoly} we obtain that its characteristic polynomial is
\[
4\bigl[(\frac{1}{2}-x)-\frac{1}{2}t+xt^2\bigr] + 2\bigl[(\frac{1}{2}-y)-\frac{1}{2}t+yt^2\bigr] = 2-3t+t^2.
\]
For the characteristic polynomial of the whole arrangement we obtain
\[
24\bigl(-\frac{1}{4}+\frac{11}{24}t-\frac{1}{4}t^2+\frac{1}{24}t^3\bigr)  = -6+11t-6t^2+t^3.
\]
Both calculations agree with \eqref{eq:braid-cpoly}.

\bibliographystyle{amsalpha}
\bibliography{char-els}

\end{document}